\documentclass[12pt,authoryear]{elsarticle}
\journal{}

\makeatletter
\def\ps@pprintTitle{%
 \let\@oddhead\@empty
 \let\@evenhead\@empty
 \def\@oddfoot{\hfill\thepage}%
 \let\@evenfoot\@oddfoot}
\makeatother

\usepackage{amsmath,amssymb,amsfonts,amsthm,graphicx}
\usepackage[bookmarksnumbered,colorlinks=true]{hyperref}
\usepackage[labelfont=bf]{caption}

\providecommand{\doi}[1]{\href{https://doi.org/#1}{DOI:#1}}
\usepackage{xurl} 
\renewcommand{\doi}[1]{%
 \href{https://doi.org/#1}{\nolinkurl{DOI:#1}}%
}

\usepackage{mathtools} 
\usepackage{geometry} 
\numberwithin{equation}{section}
\numberwithin{table}{section}
\numberwithin{figure}{section}

\theoremstyle{plain}
\newtheorem{theorem}{Theorem}[section]

\newtheorem{lemma}[theorem]{Lemma}

\theoremstyle{definition}

\newcommand{\Z}{\mathbb{Z}}
\newcommand{\N}{\mathbb{N}}
\newcommand{\R}{\mathbb{R}}
\newcommand{\PP}{\mathsf{P}} 
\newcommand{\EE}{\mathsf{E}} 
\newcommand{\leqdef}{\vcentcolon=}

\allowdisplaybreaks

\begin{document}

\begin{frontmatter}

\title{On a conjecture of Lamkin and Tkocz: Log-convexity of moments of Bernoulli sample means}

\author[a1]{Fr\'ed\'eric Ouimet\corref{mycorrespondingauthor}}

\address[a1]{Universit\'e du Qu\'ebec \`a Trois-Rivi\`eres, Trois-Rivi\`eres, QC, Canada}

\cortext[mycorrespondingauthor]{Corresponding author. Email address: frederic.ouimet2@uqtr.ca}

\begin{abstract}
Let $X_1,X_2,\ldots$ be independent $\operatorname{Bernoulli}(\theta)$ random variables, and let $\bar X_n = n^{-1}(X_1 + \cdots + X_n)$. We prove that, for every real $p \geq 1$, the sequence $\{\EE(\bar X_n^p)\}_{n \geq 1}$ is log-convex. This proves the Bernoulli case of a conjecture of Lamkin and Tkocz [Canad. Math. Bull., 65(2):271–278, 2022]. The proof conditions on the total number of successes among $2n$ trials and reduces the desired inequality to a convex-order comparison between two normalized quadratic functions of hypergeometric random variables. All the log-convexity inequalities are strict for $p > 1$ and $0 < \theta < 1$.
\end{abstract}

\begin{keyword} 
Bernoulli distribution, convex order, hypergeometric distribution, log-convexity, moments of averages
\MSC[2020]{Primary: 60E15; Secondary: 05A20, 26D15, 60C05}
\end{keyword}

\end{frontmatter}

\section{Introduction}\label{sec:intro}

Let $X_1,X_2,\ldots$ be independent copies of a nonnegative random variable, let $\bar X_n = n^{-1}(X_1 + \cdots + X_n)$, and, for $p \geq 1$ such that $\EE(X_1^p) < \infty$, consider the sequence $b_n(p) = \EE(\bar X_n^p)$. The sample means decrease in convex order. Indeed, $\bar X_{n + 1}$ is the average of the $n + 1$ leave-one-out sample means, each of which has the same distribution as $\bar X_n$, so Jensen's inequality gives $\EE[\varphi(\bar X_{n + 1})] \leq \EE[\varphi(\bar X_n)]$ for every convex function $\varphi$ for which the expectations exist. In particular, $b_n(p)$ is nonincreasing. A finer three-term property is log-convexity, namely $b_n(p)^2 \leq b_{n - 1}(p)b_{n + 1}(p)$ for every $n \geq 2$. For general background on log-convex sequences, see, for example, \citet{LiuWang2007}, and for convex order, see \citet[Chapter~3]{ShakedShanthikumar2007}.

\citet{LamkinTkocz2022} initiated a study of the logarithmic behavior of moments of averages. They proved log-convexity for $p < 0$ when the summands are positive and log-concavity for $0 < p < 1$ when the summands are nonnegative. For $p > 1$, they conjectured log-convexity for every nonnegative distribution. They also showed that the conjecture holds for integer $p \geq 1$ once $n \geq p^2$, and they singled out the Bernoulli distribution as a natural unresolved case.

More recently, \citet{JiaoLi2026} proved log-convexity for every real $p \geq 1$ in the Gamma case and for every integer $p \geq 1$ in the Poisson case.

After the first version of this paper appeared, \citet{GoncalvesMadridRamos2026} proved the Lamkin--Tkocz conjecture in full generality; the Bernoulli-specific proof presented here was obtained independently.

The purpose of this paper is to present a Bernoulli-specific proof of the conjecture. We prove log-convexity for every real $p \geq 1$ and every Bernoulli parameter $\theta \in [0,1]$, with strict inequalities when $p > 1$ and $0 < \theta < 1$. After representing the products $b_{n - 1}(p)b_{n + 1}(p)$ and $b_n(p)^2$ using two different splittings of the same collection of $2n$ Bernoulli trials, we condition on the total number of successes. The resulting problem is a convex-order comparison between two normalized quadratic functions of hypergeometric random variables. The proof uses symmetry, an interlacing of the two discrete supports, and two exact telescoping identities for the relevant probability masses.

The rest of the paper is organized as follows. Section~\ref{sec:definitions} fixes the notation and recalls the stop-loss characterization of convex order. The main result is stated in Section~\ref{sec:results} along with the hypergeometric comparison. The proofs are collected in Section~\ref{sec:proofs}.

\section{Definitions and notation}\label{sec:definitions}

Throughout, $\N = \{1,2,\ldots\}$, $\Z$ denotes the set of integers, and $x_+ = \max\{x,0\}$ for $x \in \R$. For integrable real-valued random variables $U$ and $V$, we write $U \succeq_{\mathrm{cx}} V$ if
\[
\EE[\varphi(U)] \geq \EE[\varphi(V)]
\]
for every convex function $\varphi$ for which both expectations exist. We then say that $U$ dominates $V$ in convex order. This implies $\EE(U) = \EE(V)$ because both $x\mapsto x$ and $x\mapsto -x$ are convex. Conversely, when $U$ and $V$ have the same mean, the stop-loss characterization states that
\begin{equation}\label{eq:definitions.stop.loss}
U \succeq_{\mathrm{cx}} V \quad \Longleftrightarrow \quad \EE[(z - U)_+] \geq \EE[(z - V)_+] \quad \text{for every } z \in \R;
\end{equation}
see, for example, \citet[Eq.~(3.A.6), Eq.~(3.A.8) and Theorem~3.A.1(b)]{ShakedShanthikumar2007}.

We write $A \sim \operatorname{Hypergeom}(N,K,m)$ when $A$ is the number of marked objects in a simple random sample of size $m$ drawn without replacement from a population of $N$ objects containing $K$ marked objects. Thus
\begin{equation}\label{eq:definitions.hypergeom.pmf}
\PP(A = a) = \frac{\binom{K}{a}\binom{N - K}{m - a}}{\binom{N}{m}}, \qquad a \in \Z,
\end{equation}
using the convention $\binom{u}{v} = 0$ whenever $u$ is a nonnegative integer and $v \notin \{0,\ldots,u\}$.

\section{Results}\label{sec:results}

The theorem below is our main result. It proves the Bernoulli case of Conjecture~1 of \citet{LamkinTkocz2022}, which was noted to be an open problem on page~273 of their paper.

\begin{theorem}\label{thm:Bernoulli}
Let $p \geq 1$ and $\theta \in [0,1]$, and let $X_1, X_2, \ldots$ be i.i.d.\ $\operatorname{Bernoulli}(\theta)$ random variables. Define the sample mean by $\bar X_n = n^{-1}(X_1 + \cdots + X_n)$. For $n \in \N$, define
\begin{equation}\label{eq:Bernoulli.moments}
b_n(p,\theta) \leqdef \EE(\bar X_n^p) = \sum_{k = 0}^n \binom{n}{k}\left(\frac{k}{n}\right)^p\theta^k(1 - \theta)^{n - k}.
\end{equation}
Then the sequence $\{b_n(p,\theta)\}_{n \geq 1}$ is log-convex. Equivalently,
\begin{equation}\label{eq:main.log.convexity}
b_n(p,\theta)^2 \leq b_{n - 1}(p,\theta)b_{n + 1}(p,\theta), \qquad n \geq 2.
\end{equation}
If $p > 1$ and $0 < \theta < 1$, then every inequality in \eqref{eq:main.log.convexity} is strict.
\end{theorem}

The proof of Theorem~\ref{thm:Bernoulli} is based on the following hypergeometric convex-order comparison.

\begin{lemma}\label{lem:hypergeom.cx}
Fix $n \geq 2$ and $t \in \{0,\ldots,2n\}$. For $m \in \{n - 1,n\}$, let $A_m \sim \operatorname{Hypergeom}(2n,t,m)$ and define
\[
Z_m \leqdef \frac{A_m(t - A_m)}{m(2n - m)}.
\]
Then
\begin{equation}\label{eq:hypergeom.convex.order}
Z_{n - 1} \succeq_{\mathrm{cx}} Z_n.
\end{equation}
\end{lemma}

\section{Proofs}\label{sec:proofs}

\subsection{Proof of Lemma~\ref{lem:hypergeom.cx}}

The claim is immediate for $t \in \{0,1\}$, since both random variables in \eqref{eq:hypergeom.convex.order} are then identically zero. Hence assume $t \geq 2$, and set
\[
h(a) = a(t - a), \qquad \lambda = \frac{n^2}{n^2 - 1}.
\]
Let $\pi_a = \PP(A_n = a)$ and $r_a = \PP(A_{n - 1} = a)$ for $a \in \Z$. Define the symmetrized mass function
\begin{equation}\label{eq:rho.averaging}
\rho_a = \frac{r_a + r_{t - a}}{2}, \qquad a \in \Z.
\end{equation}
Since $h(a) = h(t - a)$, a change of index gives, for every function $f$ for which the sums exist,
\begin{equation}\label{eq:symmetry}
\sum_{a \in \Z} r_a f\{h(a)\} = \frac{1}{2}\sum_{a \in \Z}\{r_a + r_{t - a}\}f\{h(a)\} = \sum_{a \in \Z}\rho_a f\{h(a)\}.
\end{equation}
Thus, when evaluating functions of $h(A_{n - 1})$, the mass function $r$ may be replaced by its symmetrization $\rho$. Moreover, the hypergeometric formula \eqref{eq:definitions.hypergeom.pmf} and the symmetry of binomial coefficients imply $\pi_a = \pi_{t - a}$ and $\rho_a = \rho_{t - a}$.

We first verify that the two variables to be compared, $\lambda h(A_{n - 1})$ and $h(A_n)$, have the same mean. The quantity $A_m(t - A_m)$ counts ordered pairs of distinct marked objects for which the first object is selected and the second is not selected. For each of the $t(t - 1)$ ordered pairs, this event has probability
\[
\frac{m}{2n}\frac{2n - m}{2n - 1}.
\]
Consequently,
\begin{equation}\label{eq:hypergeom.quadratic.mean}
\EE[A_m(t - A_m)] = \frac{t(t - 1)m(2n - m)}{2n(2n - 1)}.
\end{equation}
Taking $m = n - 1$ and $m = n$ in \eqref{eq:hypergeom.quadratic.mean}, and recalling that symmetrization does not change the law of $h(A_{n - 1})$, yields
\begin{equation}\label{eq:equal.means}
\begin{aligned}
\lambda\sum_{a \in \Z}\rho_a h(a)
&= \frac{n^2}{n^2 - 1}\EE[A_{n - 1}(t - A_{n - 1})] \\
&= \frac{n^2}{n^2 - 1}\frac{t(t - 1)(n - 1)(n + 1)}{2n(2n - 1)}
= \frac{t(t - 1)n^2}{2n(2n - 1)}
= \EE[A_n(t - A_n)] = \sum_{a \in \Z}\pi_a h(a).
\end{aligned}
\end{equation}

We now derive two telescoping identities for the masses $\pi_a$ and $\rho_a$. Set
\[
a_0 = \max\{0,t - n - 1\}.
\]
On the lower half of the support, both $\pi_a$ and $\rho_a$ vanish for $a < a_0$, and
\[
\pi_a = \frac{\binom{t}{a}\binom{2n - t}{n - a}}{\binom{2n}{n}}, \qquad a \in \Z.
\]
Define
\begin{equation}\label{eq:G.definition}
G_a = \frac{t(t - 2a - 1)}{2n(2n - t + 1)\binom{2n}{n}}\binom{t - 1}{a}\binom{2n - t + 1}{n - a}
\end{equation}
and
\begin{equation}\label{eq:K.definition}
K_a = \frac{1}{2}\pi_a h(a)(t - 2a + 1).
\end{equation}
We shall show that, for $a_0 \leq a < \lfloor t/2\rfloor$,
\begin{align}
\rho_a - \pi_a &= G_a - G_{a - 1}, \label{eq:first.telescope}\\
n^2(t - 2a - 1)G_a - h(a)\pi_a &= K_{a + 1} - K_a. \label{eq:second.telescope}
\end{align}

To prove \eqref{eq:first.telescope}, first suppose that $\pi_a > 0$. Using the hypergeometric mass function \eqref{eq:definitions.hypergeom.pmf} and adjacent-binomial ratios, we obtain
\begin{align}
\frac{r_a}{\pi_a}
&= \frac{\PP(A_{n - 1} = a)}{\PP(A_n = a)}
= \frac{\binom{t}{a}\binom{2n - t}{n - 1 - a}/\binom{2n}{n - 1}}{\binom{t}{a}\binom{2n - t}{n - a}/\binom{2n}{n}}
= \frac{n + 1}{n} \times \frac{n - a}{n - t + a + 1}, \label{eq:ra.pi.ratio}\\
\frac{r_{t - a}}{\pi_a}
&= \frac{\PP(A_{n - 1} = t - a)}{\PP(A_n = a)}
= \frac{\binom{t}{t - a}\binom{2n - t}{n - 1 - t + a}/\binom{2n}{n - 1}}{\binom{t}{a}\binom{2n - t}{n - a}/\binom{2n}{n}}
= \frac{n + 1}{n} \times \frac{n - t + a}{n - a + 1}. \label{eq:rta.pi.ratio}
\end{align}
Averaging \eqref{eq:ra.pi.ratio} and \eqref{eq:rta.pi.ratio} using \eqref{eq:rho.averaging} gives
\begin{equation}\label{eq:rho.pi.ratio}
\frac{\rho_a}{\pi_a} = \frac{n + 1}{2n}\left\{\frac{n - a}{n - t + a + 1} + \frac{n - t + a}{n - a + 1}\right\}.
\end{equation}
With $s = t - 2a$, the definition of $G_a$ similarly yields
\begin{equation}\label{eq:G.pi.ratios}
\begin{aligned}
\frac{G_a}{\pi_a}
&= \frac{t(s - 1)}{2n(2n - t + 1)}\frac{\binom{t - 1}{a}}{\binom{t}{a}}\frac{\binom{2n - t + 1}{n - a}}{\binom{2n - t}{n - a}} \\
&= \frac{t(s - 1)}{2n(2n - t + 1)}\frac{t - a}{t}\frac{2n - t + 1}{n - t + a + 1} = \frac{(t - a)(s - 1)}{2n(n - t + a + 1)}, \\
\frac{G_{a - 1}}{\pi_a}
&= \frac{t(s + 1)}{2n(2n - t + 1)}\frac{\binom{t - 1}{a - 1}}{\binom{t}{a}}\frac{\binom{2n - t + 1}{n - a + 1}}{\binom{2n - t}{n - a}} \\
&= \frac{t(s + 1)}{2n(2n - t + 1)}\frac{a}{t}\frac{2n - t + 1}{n - a + 1} = \frac{a(s + 1)}{2n(n - a + 1)}.
\end{aligned}
\end{equation}
Subtracting $1$ from \eqref{eq:rho.pi.ratio}, putting the two terms over a common denominator, and using \eqref{eq:G.pi.ratios}, we obtain
\begin{equation}\label{eq:first.telescope.ratio.form}
\begin{aligned}
\frac{\rho_a - \pi_a}{\pi_a}
&= \frac{n + 1}{2n}\left\{\frac{n - a}{n - t + a + 1} + \frac{n - t + a}{n - a + 1}\right\} - 1 \\
&= \frac{1}{2n}\left\{\frac{(n + 1)(n - a)}{n - t + a + 1} + \frac{(n + 1)(n - t + a)}{n - a + 1} - 2n\right\} \\
&= \frac{1}{2n}\frac{(t - a)(s - 1)(n - a + 1) - a(s + 1)(n - t + a + 1)}{(n - t + a + 1)(n - a + 1)} \\
&= \frac{1}{2n}\left\{\frac{(t - a)(s - 1)}{n - t + a + 1} - \frac{a(s + 1)}{n - a + 1}\right\} = \frac{G_a - G_{a - 1}}{\pi_a}.
\end{aligned}
\end{equation}
This proves \eqref{eq:first.telescope} whenever $\pi_a > 0$. If $\pi_a = 0$ for some $a_0 \leq a < \lfloor t/2\rfloor$, then necessarily $t \geq n + 1$ and $a = a_0 = t - n - 1$. At this point $r_{t - a_0} = 0$, and direct substitution gives
\[
\rho_{a_0} = G_{a_0} = \frac{n + 1}{2n\binom{2n}{n}}\binom{t}{a_0}, \qquad G_{a_0 - 1} = 0.
\]
Thus \eqref{eq:first.telescope} also holds at the boundary.

For \eqref{eq:second.telescope}, the adjacent-mass ratio is
\begin{equation}\label{eq:adjacent.pi.ratio}
\frac{\pi_{a + 1}}{\pi_a} = \frac{(t - a)(n - a)}{(a + 1)(n - t + a + 1)}
\end{equation}
whenever $\pi_a > 0$. Writing again $s = t - 2a$, we have $h(a + 1) = (a + 1)(t - a - 1)$ and $t - 2(a + 1) + 1 = s - 1$. The algebraic identity
\[
(n - a)(t - a - 1)(s - 1) - a(s + 1)(n - t + a + 1) = n(s - 1)^2 - 2a(n - t + a + 1),
\]
where $t = 2a + s$, will also be used below. Hence equations \eqref{eq:K.definition} and \eqref{eq:adjacent.pi.ratio} imply
\begin{align*}
\frac{K_{a + 1} - K_a}{\pi_a}
&= \frac{t - a}{2(n - t + a + 1)}\left\{(n - a)(t - a - 1)(s - 1) - a(s + 1)(n - t + a + 1)\right\} \\
&= \frac{n(t - a)(s - 1)^2}{2(n - t + a + 1)} - a(t - a) \\
&= \frac{n^2(s - 1)G_a}{\pi_a} - h(a),
\end{align*}
where the last equality follows from the first ratio in \eqref{eq:G.pi.ratios}. This proves \eqref{eq:second.telescope} whenever $\pi_a > 0$. At the possible boundary point $a = a_0$ identified above, one has $K_{a_0} = 0$ and
\[
K_{a_0 + 1} = \frac{n(n + 1)(t - 2a_0 - 1)}{2\binom{2n}{n}}\binom{t}{a_0} = n^2(t - 2a_0 - 1)G_{a_0},
\]
so \eqref{eq:second.telescope} holds there as well. Notice also that
\begin{equation}\label{eq:boundary.and.positivity}
G_{a_0 - 1} = 0, \qquad K_{a_0} = 0, \qquad G_a \geq 0 \quad \text{for } a_0 \leq a < \lfloor t/2\rfloor.
\end{equation}
Indeed, if $a_0 = 0$, then $h(a_0) = 0$, while if $a_0 > 0$, then $\pi_{a_0} = 0$; this proves the middle equality in \eqref{eq:boundary.and.positivity}. If $a_0 = 0$, then $\binom{t - 1}{a_0 - 1} = \binom{t - 1}{-1} = 0$. If $a_0 > 0$, then $a_0 = t - n - 1$ and
\[
n - (a_0 - 1) = 2n - t + 2 > 2n - t + 1,
\]
so the second binomial coefficient in the definition of $G_{a_0 - 1}$ vanishes. This proves the first equality in \eqref{eq:boundary.and.positivity}. The final assertion follows from \eqref{eq:G.definition} and $t - 2a - 1 > 0$ for $a < \lfloor t/2\rfloor$.

The next ingredient is an interlacing of the support points. Put
\[
y_a = h(a), \qquad x_a = \lambda y_a.
\]
We claim that
\begin{equation}\label{eq:interlace}
y_a \leq x_a \leq y_{a + 1}, \qquad a_0 \leq a < \lfloor t/2\rfloor.
\end{equation}
The first inequality follows from $\lambda > 1$ and $y_a \geq 0$. For the second inequality in \eqref{eq:interlace}, note that $y_{a + 1} - y_a = s - 1$, $y_a = (t^2 - s^2)/4$, $t \leq 2n$, and $s\geq 2$, so we have
\[
\begin{aligned}
(n^2 - 1)(y_{a + 1} - x_a)
&= (n^2 - 1)(y_{a + 1} - y_a) - y_a \\
&= (n^2 - 1)(s - 1) - \frac{t^2 - s^2}{4} \\
&\geq (n^2 - 1)(s - 1) - n^2 + \frac{s^2}{4}
= (s - 2)\left(n^2 + \frac{s - 2}{4}\right) \geq 0.
\end{aligned}
\]
This proves \eqref{eq:interlace}.

It remains to verify the stop-loss inequalities in \eqref{eq:definitions.stop.loss}. Define
\begin{equation}\label{eq:D.definition}
D(z) = \sum_{a \in \Z}\rho_a(z - \lambda y_a)_+ - \sum_{a \in \Z}\pi_a(z - y_a)_+, \qquad z \in \R.
\end{equation}
For $a_0 \leq c < \lfloor t/2\rfloor$, equation \eqref{eq:first.telescope} and summation by parts give
\begin{align}
\sum_{a = a_0}^c(\rho_a - \pi_a) &= G_c, \label{eq:partial.mass.sum}\\
\sum_{a = a_0}^c y_a(\rho_a - \pi_a) &= y_cG_c - \sum_{a = a_0}^{c - 1}(y_{a + 1} - y_a)G_a = y_cG_c - \sum_{a = a_0}^{c - 1}(t - 2a - 1)G_a. \label{eq:partial.weighted.sum}
\end{align}
Because $c < \lfloor t/2\rfloor$, the masses contributing to \eqref{eq:D.definition} at $z = x_c$ occur in symmetric pairs. Thus every active lower-half contribution has a distinct reflected counterpart; this accounts for the factor $2$ in \eqref{eq:D.at.xc.preliminary} and \eqref{eq:D.slope}. Using \eqref{eq:partial.mass.sum} and \eqref{eq:partial.weighted.sum}, we obtain
\begin{equation}\label{eq:D.at.xc.preliminary}
\begin{aligned}
\frac{D(x_c)}{2}
&= \lambda y_c\sum_{a = a_0}^c(\rho_a - \pi_a) - \lambda\sum_{a = a_0}^c y_a\rho_a + \sum_{a = a_0}^c y_a\pi_a \\
&= \lambda\sum_{a = a_0}^{c - 1}(t - 2a - 1)G_a - (\lambda - 1)\sum_{a = a_0}^c y_a\pi_a.
\end{aligned}
\end{equation}
Multiplying \eqref{eq:D.at.xc.preliminary} by $n^2 - 1$ and then using \eqref{eq:second.telescope} and $K_{a_0} = 0$ gives
\[
\frac{n^2 - 1}{2}D(x_c)
= \sum_{a = a_0}^{c - 1}\{n^2(t - 2a - 1)G_a - y_a\pi_a\} - y_c\pi_c
= K_c - y_c\pi_c = \frac{1}{2}\pi_cy_c(t - 2c - 1).
\]
Therefore,
\begin{equation}\label{eq:D.at.xc.nonnegative}
D(x_c) = \frac{\pi_cy_c(t - 2c - 1)}{n^2 - 1} \geq 0, \qquad a_0 \leq c < \lfloor t/2\rfloor.
\end{equation}

The interlacing in \eqref{eq:interlace} orders the candidate breakpoints as
\[
y_{a_0} \leq x_{a_0} \leq y_{a_0 + 1} \leq x_{a_0 + 1} \leq \cdots \leq y_{\lfloor t/2\rfloor} \leq x_{\lfloor t/2\rfloor}.
\]
The initial interval $[y_{a_0},x_{a_0}]$ causes no difficulty. If $a_0 = 0$, then $y_{a_0} = x_{a_0} = 0$. If $a_0 > 0$, then $\pi_{a_0} = 0$, and no term in \eqref{eq:D.definition} is active before $x_{a_0}$, so $D$ vanishes throughout this interval. On each open interval $(x_c,y_{c + 1})$, symmetry and \eqref{eq:partial.mass.sum} show that the slope of $D$ is
\begin{equation}\label{eq:D.slope}
D'(z) = 2\sum_{a = a_0}^c(\rho_a - \pi_a) = 2G_c \geq 0.
\end{equation}
Hence $D(y_{c + 1}) \geq D(x_c) \geq 0$. On every interval $[y_{c + 1},x_{c + 1}]$, the function $D$ is affine, and both endpoint values are nonnegative by the preceding inequality and \eqref{eq:D.at.xc.nonnegative}, except at the final endpoint $x_{\lfloor t/2\rfloor}$, which is handled next. Since $x_{\lfloor t/2\rfloor}$ is at least as large as every support point and the means agree by \eqref{eq:equal.means},
\begin{equation}\label{eq:D.at.final.point}
\begin{aligned}
D(x_{\lfloor t/2\rfloor})
&= x_{\lfloor t/2\rfloor}\left(\sum_{a \in \Z}\rho_a - \sum_{a \in \Z}\pi_a\right) - \lambda\sum_{a \in \Z}\rho_a y_a + \sum_{a \in \Z}\pi_a y_a \\
&= 0.
\end{aligned}
\end{equation}
Also, $D(z) = 0$ below the smallest support point and for $z \geq x_{\lfloor t/2\rfloor}$. The function $D$ is continuous and affine between consecutive breakpoints. Equations \eqref{eq:D.at.xc.nonnegative}, \eqref{eq:D.slope}, and \eqref{eq:D.at.final.point} therefore imply
\begin{equation}\label{eq:D.nonnegative}
D(z) \geq 0, \qquad z \in \R.
\end{equation}

By \eqref{eq:equal.means}, \eqref{eq:D.nonnegative}, and the stop-loss characterization \eqref{eq:definitions.stop.loss},
\[
\lambda h(A_{n - 1}) \succeq_{\mathrm{cx}} h(A_n).
\]
Here, by \eqref{eq:symmetry}, the distribution of $h(A_{n - 1})$ may equivalently be represented using the symmetrized mass function $\rho$. Finally,
\[
Z_{n - 1} = \frac{\lambda h(A_{n - 1})}{n^2}, \qquad Z_n = \frac{h(A_n)}{n^2}.
\]
Convex order is preserved under multiplication by a positive constant, so \eqref{eq:hypergeom.convex.order} follows.

\subsection{Proof of Theorem~\ref{thm:Bernoulli}}

Fix $n \geq 2$. The cases $\theta \in \{0,1\}$ are immediate, so assume $0 < \theta < 1$, and let $Y_1,\ldots,Y_{2n}$ be independent $\operatorname{Bernoulli}(\theta)$ random variables. Set
\[
T = \sum_{i = 1}^{2n}Y_i, \qquad A_m = \sum_{i = 1}^mY_i, \qquad m \in \{n - 1,n\}.
\]
Unconditionally, the variables $A_m$ and $T - A_m$ are independent binomial random variables based on disjoint groups of sizes $m$ and $2n - m$. Therefore, with the notation in \eqref{eq:Bernoulli.moments},
\begin{equation}\label{eq:block.product.identity}
b_m(p,\theta)b_{2n - m}(p,\theta)
= \EE\left[\left(\frac{A_m}{m}\right)^p\left(\frac{T - A_m}{2n - m}\right)^p\right]
= \EE\left[\left\{\frac{A_m(T - A_m)}{m(2n - m)}\right\}^p\right].
\end{equation}
Conditionally on $T = t$, exchangeability gives
\[
A_m \mid (T = t) \sim \operatorname{Hypergeom}(2n,t,m).
\]
Lemma~\ref{lem:hypergeom.cx} and the convexity of $x \mapsto x^p$ on $[0,\infty)$ now imply, for every $t \in \{0,\ldots,2n\}$,
\begin{equation}\label{eq:conditional.moment.comparison}
\EE\left[\left.\left\{\frac{A_{n - 1}(t - A_{n - 1})}{n^2 - 1}\right\}^p\right|T = t\right]
\geq \EE\left[\left.\left\{\frac{A_n(t - A_n)}{n^2}\right\}^p\right|T = t\right].
\end{equation}
Averaging \eqref{eq:conditional.moment.comparison} over $T$ and applying \eqref{eq:block.product.identity} with $m = n - 1$ and $m = n$ yields
\[
b_{n - 1}(p,\theta)b_{n + 1}(p,\theta) \geq b_n(p,\theta)^2,
\]
which proves log-convexity.

It remains to verify strictness. Assume $p > 1$ and $0 < \theta < 1$. Then $\PP(T = 2) > 0$. Given $T = 2$, the quantity $A_m(2 - A_m)$ equals $1$ when $A_m = 1$ and equals $0$ otherwise. Moreover,
\[
\PP(A_m = 1 \mid T = 2) = \frac{2\binom{2n - 2}{m - 1}}{\binom{2n}{m}} = \frac{m(2n - m)}{n(2n - 1)}.
\]
Consequently,
\begin{equation}\label{eq:T2.moment}
\EE\left[\left.\left\{\frac{A_m(2 - A_m)}{m(2n - m)}\right\}^p\right|T = 2\right] = \frac{\{m(2n - m)\}^{1 - p}}{n(2n - 1)}.
\end{equation}
Since $1 - p < 0$ and $n^2 - 1 < n^2$, the right-hand side of \eqref{eq:T2.moment} is strictly larger for $m = n - 1$ than for $m = n$. Thus \eqref{eq:conditional.moment.comparison} is strict on an event of positive probability, and the averaged inequality is strict. For completeness, when $p = 1$, one has $b_n(1,\theta) = \theta$ for every $n$, and when $\theta \in \{0,1\}$, the sequence is constant; these are precisely the elementary equality cases covered by the theorem.

\setlength{\bibsep}{0pt plus 0ex}


\section*{Statements and declarations}
\addcontentsline{toc}{section}{Statements and Declarations}

\subsection*{Competing interests}
\addcontentsline{toc}{section}{Competing interests}

The author has no relevant financial or non-financial interests to disclose.

\subsection*{Funding}
\addcontentsline{toc}{section}{Funding}

The author is supported by the Natural Sciences and Engineering Research Council of Canada (Discovery Grant RGPIN-2026-04471, Discovery Launch Supplement DGECR-2026-00449).

\subsection*{Data availability statement}
\addcontentsline{toc}{section}{Data availability statement}

No data were used or generated in this study.

\end{document}